\documentclass{article}

\usepackage{graphicx}
\usepackage{amssymb,amsmath}

\newcommand{\de}{{\rm d}}

\title{Nonparametric estimation of distribution and density functions in 
  presence of missing data: an IFS approach}

\author{Stefano M. Iacus\footnote{Department of Economics, Via Conservatorio 7, I-20122 Milan - Italy, email: stefano.iacus@unimi.it} \and Davide La Torre}

\DeclareGraphicsExtensions{.eps}

\begin{document}

\maketitle

\begin{abstract}
% Text of abstract
In this paper we consider a class of nonparametric estimators of a distribution function $F$, with compact support, based on the theory of IFSs. The estimator
of $F$ is tought as the fixed point of a
contractive operator $T$ defined in terms of a vector of parameters $p$ and a family of affine maps $\mathcal W$
which can be both depend of the sample $(X_1, X_2, \ldots, X_n)$.
Given $\mathcal W$, the problem consists in finding
a vector $p$ such that the fixed point of $T$ is ``sufficiently
near'' to $F$. It turns out that this is a quadratic constrained optimization problem that we propose to  solve by  penalization techniques.  If $F$ has a density $f$, we can also provide an estimator of $f$ based on Fourier techniques. IFS estimators for $F$ are asymptotically equivalent to the empirical distribution function (e.d.f.) estimator. We will study relative efficiency of the IFS estimators with respect to the e.d.f. for small samples via Monte Carlo approach.

For well behaved distribution functions $F$ and for a particular family of so-called wavelet maps
the IFS estimators can be dramatically better than the e.d.f. (or the kernel estimator for density estimation) in presence of missing
data, i.e. when it is only possibile to observe data on subsets of the whole support of $F$.

This research has also produced a free package for the R statistical environment which is ready to be used in applications.
\end{abstract}

\textbf{key words:} iterated function systems, distribution function estimation,
nonparametric estimation, missing data, density estimation.\\

% main text
\section{Introduction}
Let $X_1, X_2, \ldots, X_n$ be an i.i.d.  sample drawn from a
random variable $X$ with unknown distribution function
$F$ with compact support $[\alpha,\beta]$. The empirical
distribution function (e.d.f.) $$ \hat F_n(x) =
\frac{1}{n}\sum_{i=1}^n \chi(X_i\leq x) $$ is one commonly used
estimator of the unknown distribution function $F$ (here $\chi$ is the indicator function). The e.d.f. has
an impressive set of good statistical properties such as it is
first order efficient in the minimax sense (see \cite{dkw}, \cite{beran}, \cite{levit},
\cite{millar}, and \cite{gilevit}). More or less recently, other
second order efficient estimators have been proposed in the
literature for special classes of distribution functions $F$.
Golubev and Levit (1996a, b) and \cite{efro} are two of such
examples. It is rather curious that a step-wise function can be
such a good estimator and, in fact, \cite{efro} shows that,
for the class of analytic functions, for small sample sizes, the
e.d.f. is not the best estimator. In this paper we study the
properties of a new class of distribution function estimators
based on iterated function systems (IFSs) introduced by the
authors in a previous work \cite{iaclat}. IFSs have
been introduced in \cite{hutch} and \cite{bd}. The main idea on which this method is based consists of
thinking the estimation of $F$ as the fixed point of a contraction $T$
on a complete metric space. The operator $T$ is defined in terms
of a family of affine  maps $\mathcal W$ and a vector of parameters $p$. For a given family $\mathcal W$, $T$ depends only on the choice $p$. The idea,
known as {\it inverse approach} (see Section \ref{sec2})  is to
determine $p$ by solving a constrained quadratic optimization
problem built in terms of sample moments. In this paper
this optimization problem is solved by a penalization method.
The nature of affine maps allow to derive easily the Fourier transform of $F$ and, when available, an explicit formula
for the density of $F$ via anti Fourier transform.
In this way, given $\mathcal W$ and $p$ we have at the same time estimators for the distribution, characteristic and density functions.

The paper is organized as follows. In Section \ref{sec2}  the
inverse approach is presented and a penalization method is proposed in
order to solve a quadratic optimization problem. We also discuss
the choice of the family of maps $\mathcal W$. In Section
\ref{sec3} numerical results and comparisons with classical
estimators are shown for small samples via Monte Carlo Analysis. 

Finally we show an application of these estimators when the empirical distribution function (or the kernel density estimator for the density) cannot be applied. We will consider situations of missing data when, for example, the data can only be observed on some windows of the support of $F$. This can be the case of directional data analysis when, for some reason, instruments are not able for technical or physical reason to collect data in same range of angles say $A$ and $B$, $A, B\subseteq [0,2\pi]$. For $x$ in $A$ or $B$ the e.d.f. will be constant and, at the same time, the kernel density estimator will estimate a plurimodal distribution for these data.
In this case we will show examples in which the IFS estimator does it job incredibly well.

Tables and figures can be found at the end of the paper after the references.

\section{An IFS estimator}\label{sec2}
The theory of distribution function approximation via IFSs we will use to derive estimators is due to  \cite{fv95}. Results from this section, apart from were explicitly mentioned, are from the cited authors. Let $\mathcal M(X)$ be the set of probability measures
on $\mathcal B(X)$, the $\sigma$-algebra of Borel subsets of $X$
where $(X,d)$ is a compact metric space (in our case will be
$X=[\alpha,\beta]$ and $d$ the Euclidean metric.)

In the IFSs literature the following {\sl Hutchinson} metric plays
a crucial role $$ d_H(\mu,\nu) = \sup_{f\in {\rm Lip}(X)} \left\{
\int_X f \de \mu - \int_X f \de \nu \right\}, \quad \mu,\nu \in
\mathcal M(X) $$ where $$ {\rm Lip}(X) = \{ f: X\to \mathbb R,
|f(x)-f(y)|\leq d(x,y), x,y \in X\} $$ thus $(\mathcal M(X),d_H)$
is a complete metric space \cite[see][]{hutch}.

We denote by $({\bf w},{\bf p})$ an {\sl $N$-maps
contractive IFS on $X$ with probabilities} or simply an {\sl
$N$-maps IFS}, that is, a set of $N$ affine contraction maps,
${\bf w} = (w_1,w_2,\ldots,w_N)$, $$w_i = a_i + b_i \, x,\quad {\rm with}\,\, |b_i|<1,\quad
b_i,a_i\in\mathbb R,\quad i=1,2,\ldots, N, $$ with
associated probabilities ${\bf p} = (p_1,p_2,\ldots,p_N)$,
$p_i\geq 0$, and $\sum_{i=1}^N p_i =1$. The IFS has a
contractivity factor defined as $$ c = \max_{1\leq i\leq N}
|b_i| <1 $$ Consider the following (usually called {\sl
Markov}) operator $M : \mathcal M(X)\to \mathcal M(X)$ defined as

\begin{equation}
M\mu = \sum_{i=1}^N p_i \mu \circ w_i^{-1},\quad \mu \in \mathcal
M(X),
\end{equation}

where $w_i^{-1}$ is the inverse function of $w_i$ and $\circ$
stands for the  composition. In Hutchinson (1981) it was shown
that $M$ is a contraction mapping on $(\mathcal M(X),d_H)$ i.e.
for all $\mu,\nu\in \mathcal M(X)$, $d_H(M\mu,M\nu)\leq c
d_H(\mu,\nu)$. Thus, there exists a unique measure
$\bar\mu\in\mathcal M(X)$, the {\sl invariant measure} of the IFS,
such that $M\bar\mu = \bar \mu$ by Banach theorem. Associated to
each measure $\mu\in \mathcal M(X)$, there exists a distribution
function $F$. In terms of it the previous operator $M$ can be
rewritten as 
$$ TF(x)=\left\{
\begin{array}{ll} 0 & \ \ \mbox{if $x\le \alpha$}\\
\\
\sum\limits_{i=1}^N p_i F(w_i^{-1}(x)) & \ \ \mbox{if $\alpha<x<\beta$}\\
\\
1 & \ \ \mbox{if $x\ge \beta$}\\

\end{array}
\right. $$

\subsection{Minimization approach}

For affine IFSs there exists a simple and useful  relation between
the moments of probability measures  on $\mathcal M(X)$. Given a
$N$-maps IFS$({\bf w},{\bf p})$ with associated Markov operator
$M$, and given a measure $\mu\in\mathcal M(X)$  then, for any
continuous function $f:X\to\mathbb R$,

\begin{equation}
\int_X f(x) \de \nu(x) = \int_X f(x) \de(M\mu)(x) = \sum_{i=1}^N
p_i \int_X (f\circ w_i)(x)\de \mu(x)\,, \label{eq1}
\end{equation}

where $\nu = M\mu$. In our case $X = [\alpha,\beta]\subset\mathbb R$  so we
readly have a relation involving the moments of $\mu$ and $\nu$.
Let
\begin{equation}
g_k = \int_X x^k \de\mu,\quad h_k = \int_X x^k \de \nu,\quad
k=0,1,2,\ldots,
\end{equation}
be the moments of the two measures, with $g_0 = h_0 = 1$. Then, by
\eqref{eq1}, with $f(x) = x^k$, we have $$ h_k = \sum_{j=0}^k
\binom{k}{j} \left\{\sum_{i=1}^N p_i b_i^j a_i^{k-j}
\right\} g_j,\quad k=1,2,\ldots,\,. $$
Set $X=[\alpha,\beta]$ and let $\mu$ and $\mu^{(j)} \in \mathcal M(X)$,
$j=1,2,\ldots$ with associated moments of any order $g_k$ and $$
g_k^{(j)}  = \int_X x^k \de \mu^{(j)}\,. $$ Then, the following
statements are equivalent (as $j\to\infty$  and $\forall k\geq
0$):
\begin{enumerate}
\item $g_k^{(j)}\to g_k$,
\item $\forall f\in {\bf C}(X)$, $\int_X f \de\mu^{(j)} \to \int_X f\de\mu\,,$ (weak* convergence),
\item $d_H(\mu^{(j)},\mu)\to0$.
\end{enumerate}

(here ${\bf C}(X)$ is the space of continuous functions on $X$).
This result gives a way to find and appropriate set of maps and
probabilities by solving the so called problem of moment matching.
With the solution in hands, given the convergence of the moments,
we also have the convergence of the measures and then the
stationary measure of $M$ approximates with given precision (in a
sense specified by the collage theorem below) the target measure
$\mu$ \cite[see][]{bd}.

Next result, called the {\sl collage} theorem is a standard
product  of the IFS theory and is a consequence of Banach theorem.

{\bf (Collage theorem) : }
Let $(Y,d_Y)$ be a complete metric space. Given an $y\in Y$,
suppose that there  exists a contractive map $f$ on $Y$ with
contractivity factor $0\leq c<1$ such that
$d_Y(y,f(y))<\varepsilon$. If $\bar y$ is the fixed point of $f$,
i.e. $f(\bar y) = \bar y$, then $d_Y(\bar y,y) <
\frac{\varepsilon}{1-c}$.

So if one wishes to approximate  a function $y$ with the fixed
point $\bar y$ of an unknown contractive map $f$, it is only
needed to solve the inverse problem of finding $f$ which minimizes
the collage distance $d_Y(y,f(y))$.

The main result in Forte and Vrscay that we will use to build one
of the IFS estimators is that the inverse problem can be reduced
to minimize a suitable quadratic form in terms of the $p_i$ given
a set of affine maps $w_i$ and the sequence of moments $g_k$ of
the target measure. Let $$\Pi^N = \left\{ {\bf p} =
(p_1,p_2,\ldots,p_N) : p_i\geq 0, \sum_{i=1}^N p_i = 1 \right\} $$
be the simplex of probabilities. Let ${\bf w} =
(w_1,w_2,\ldots,w_N)$, $N=1,2,\ldots$ be  subsets of $\mathcal W =
\{w_1,w_2,\ldots\}$ the infinite set of affine contractive maps on
$X=[\alpha,\beta]$ and let ${\bf g}$ the set of the moments of any order of
$\mu\in\mathcal M(X)$. Denote by $M$ the Markov operator of the
$N$-maps IFS $({\bf w},{\bf p})$ and by  $\nu_N = M\mu$, with
associated moment vector of any order ${\bf h}_N$. The collage
distance between the moment vector of $\mu$ and $\nu_N$ $$
\Delta({\bf p}) = ||{\bf g}-{\bf h}_N||_{\bar l^2} : \Pi^N \to
\mathbb R $$ is a continuous function and attains an absolute
minimum value $\Delta_{\min}$ on $\Pi^N$ where $$ ||{\bf
x}||_{\bar l^2} = x_0^2 + \sum_{k=1}^\infty \frac{x_k^2}{k^2}\,.$$
Moreover,
$\Delta^N_{\min} \to 0$ as $N\to\infty$.
Thus, the collage distance can be made arbitrarily small by
choosing a suitable number of maps and probabilities.

The above inverse problem can be posed as a
quadratic programming one in the following notation $$ S({\bf p})
= (\Delta({\bf p}))^2 = \sum_{k=1}^\infty \frac{(h_k-g_k)^2}{k^2}
$$ $$D(X) = \{{\bf g} = (g_0,g_1,\ldots) : g_k = \int_X x^k
\de\mu, k=0,1,\ldots, \mu\in\mathcal M(X)\} $$

Then by \eqref{eq1} there exists a linear operator $A :D(X)\to
D(X)$ associated to $M$ such that ${\bf h}_N = A {\bf g}$. In
particular 
\begin{equation}
h_k = \sum_{i=1}^N A_{ki} p_i,\quad k=1,2,\ldots
\quad\text{where} \quad
A_{ki} = \sum_{j=0}^\infty \binom{k}{j} b_i^j a_i^{k-j} g_j
\label{eqAik}
\end{equation}
Thus $$ S({\bf p}) = {\bf p}^t Q {\bf p} + {\bf B}^t {\bf p}
+C,\leqno{({\bf Q})} $$ 
$$\text{where}\quad Q=[q_{ij}],\quad   q_{ij} = \sum_{k=1}^\infty
\frac{A_{ki} A_{kj}}{k^2},\quad i,j = 1,2,\ldots, N, $$ 
\begin{equation}
B_i =
-2\sum_{k=1}^\infty \frac{g_k}{k^2} A_{ki},\quad i=1,2,\ldots, N
\quad \text{and}\quad C =\sum_{k=1}^\infty \frac{g_k^2}{k^2}\,.
\label{eqBiC}
\end{equation}
The series above are convergent as $0\leq A_{ni}\leq 1$ and the
minimum can be found by minimizing the quadratic form on the
simplex $\Pi^N$.

The estimator will then be built by substituting the moments of the target measure with the empirical moments and by truncation of the above series to a finite sum.

\subsection{Numerical solutions}
When practical cases are considered, in particular concerning
estimation, the previous series have to be truncated and this
implies that the matrix $Q$ is assured to be definite positive. Standard
numerical procedures for the minimization of constrained quadratic
optimization problems involving positive definite quadratic forms
cannot be used in this context. To solve this problem an approach
is to build the following penalized function $L_\lambda$

$$ L_\lambda({\bf p})={\bf p}^t Q {\bf p} + {\bf B}^t {\bf p}
+C+\lambda\left(1-\sum_{i=1}^N p_i\right)^2 $$

and then to study the following problem

$$ \min L_\lambda({\bf p}), \ \  0\le p_i\le 1  \leqno{(LOP)} $$

It is trivial that an optimizer ${\bf p^*}$ of (LOP) such that
$\sum_{i=1}^N p_i^*=1$ is also an optimizer for the problem

$$ \min S({\bf p}), \ \  {\bf p}\in \Pi^N  \leqno{(OP)} $$

For solving (LOP) numerically, we have used the method L-BFGS-B due to \cite{byrd} which allows to minimize a nonlinear function with box
constraints, i.e. when each variable can be given a lower and/or
upper bound. The initial value of this procedure must satisfy the
constraints. This uses a limited-memory modification of the BFGS
quasi-Newton method. The method `''BFGS''' is a quasi-Newton method
(also known as a variable metric algorithm).

\subsection{The choice of affine maps}

As we are mostly concerned with estimation, we briefly discuss the
problem of choosing the maps. In \cite{fv95} the
following two sets of wavelet-type maps are proposed. Fixed and
index $i^*\in\mathbb N$, define

$$ \gamma_{ij} = \frac{x-\alpha+(j-1)(\beta-\alpha)}{2^i}+\alpha,\quad i=1,2,\ldots,
i^* \quad j = 1,2,\ldots, 2^i $$ and $$ \eta_{ij} =
\frac{x-\alpha+(j-1)(\beta-\alpha)}{i},\quad i=2,\ldots, i^* \quad j = 2,\ldots,
i\,. $$ Then set $N = \sum\limits_{i=1}^{i^*} 2^i$ or
$N=i^*(i^*-1)/2$ respectively. To choose the maps, consider the
natural ordering of the maps $\omega_{ij}$ and operate as follows
$$ \mathcal W_1 =\{ w_1 = \gamma_{11}, w_2 = \gamma_{12}, w_3 =
\gamma_{21}, \ldots, w_6 =\gamma_{24}, w_7 = \gamma_{31}, \ldots,
w_{N}=\gamma_{i^*2^{i^*}}\} $$ and $$\mathcal W_2 =\{ w_1 =
\eta_{22}, w_2 = \eta_{32}, w_3 = \eta_{33}, w_4 =\eta_{42},
\ldots, w_6 = \eta_{44}, \ldots, w_{N}=\eta_{i^*i^*}\} $$
respectively. In \cite{iaclat} we proposed the
following quantile based maps $$\mathcal Q_1 =\{ w_i (x) =
(q_{i+1}-q_i) x + q_i, i=1,2,\ldots, N\}$$ where $q_i =
F^{-1}(u_i)$, and $0=u_1 < u_2 < \ldots < u_{N} < u_{N+1} = 1$ are
$N+1$ equally spaced points on $[0,1]$.
With these maps, it has been shown that, there is no need to use a moment matching approach. In particular, given $p_i=1/N$, the IFSs turns out to be a smoother of the e.d.f. and so it has nice small sample and asymptotic statistical properties (see cited reference) even for non compact support distribution functions $F$.
Here we will also mix the quantile information with the moment matching idea. To distinguish the two cases (fixed $p_i=1/N$ or $p$ solution of $({\bf QP})$) we will use the notation $\mathcal Q_1$ and $\mathcal Q_2$ later on.

\subsection{Fourier analysis results}
We recall, from \cite{fv98}  results that are rather straight forward to prove but also essential  to us since we will use them in density estimation and in particular in presence of missing
data. Simplicity is due to affinity of the maps. We assume that the support of the measures is $X= [0,1]$ without loss of generality.

Given a measure $\mu\in\mathcal M(X)$, the Fourier transform (FT)  $\phi : \mathbb R \to \mathbb C$, where $\mathbb C$ is the complex space, is defined by the relation
$$
\phi(t) = \int_X e^{-itx} \de \mu(x),\quad t\in\mathbb R\,,
$$
with the well known properties $\phi(0) =1$ and $|\phi(t)|\leq 1$, $\forall\, t\in\mathbb R$.
It can be shown that the space of characteristic functions ${\mathcal FT}(X)$ can be made metric and complete with an opportune metric.
Thus, given a $N$-maps affine IFS$({\bf w},{\bf p})$  it is possibile to define a new linear operator $B: {\mathcal FT}(X)\to {\mathcal FT}(X)$ whose unique fixed point
reads as
$$
\bar\phi(t) = \sum_{k=1}^N p_k e^{-i t a_k} \bar\phi(b_k t),\quad t\in\mathbb R\,.
$$
This $\bar\phi(t)$ is the FT of the fixed point of the $N$-maps affine IFS$({\bf w},{\bf p})$.
Now \cite[see e.g.][]{tarter}, suppose that the target distribution $F$ admits a density $f$. It is possible to write the density $f$ via Fourier expansion. In fact,
$$
\phi(t) =\int_0^1 f(x)  e^{-itx}  \de x = \int_0^1 e^{-itx} \de F(x)
$$
thus
$$
f(x) = \frac{1}{2\pi}\sum_{k=-\infty}^{+\infty} B_k e^{ikx} 
\quad\text{where}\quad B_k  = \phi(k)\,. 
$$

\section{Relative efficiency and estimation in presence of missing data}\label{sec3}
Suppose to have an i.i.d. sample on $n$ observations with common unknown distribution function $F$ with compact support on $[\alpha,\beta]$ which has all the moments up to order $M$. An IFS estimator of $F$ is the fixed point of the functional $TF$ where the $N$ maps are choosen in advance and the $p_i$ are the solution of the ({\bf QP}) quadratic programming problem where in the expression on $A_{ik}$, $B_i$ and $C$ we replace, in equations \eqref{eqBiC}
 and \eqref{eqAik},  the true moments $g_k$ with the sample moments $m_k$, $k=0,1,\ldots, M$ for a fixed $M$ and we consider the first $M$ terms of the series involved.

Given the solution of ({\bf QP}), we have an estimator for $F$ and an estimator for the characteristic function of $F$, say $\hat \phi$. Suppose that $F$ posseses a density $f$ then we have further a (Fourier) density estimator for $f$
$$
\begin{aligned}
\hat f(x) &=  \frac{1}{2\pi}\sum_{k=-m}^{+m} \hat B_k e^{ikx}\\
&=
\frac{1}{2\pi} + \frac{1}{\pi} \sum_{i=1}^m \biggl\{{\rm Re}(\hat B_k)\,\cos(kx) -
 {\rm Im}(\hat B_k)\,\sin(kx)\biggr\}
\end{aligned}
$$
where $\hat B_k = \hat \phi(k)$ and $m$, the number of Fourier terms, is choosen in the usual way,
i.e.
$$
\text{if} \, \left | \hat B_{m+1} \right |^2 \text{and}   \left | \hat B_{m+2} \right |^2 < \frac{2}{n+1}
\quad\text{then use the first $m$ coefficents}
$$
\cite[see again][]{tarter}.
Tables \ref{tab:a} and \ref{tab:b} show camparisons between the empirical cumulative distribution function $\hat F_n$ and the IFS estimator, say $\hat T_N$,  for some target distributions $F$, in terms of average mean square error (AMSE) and sup-norm (SUP) distance.
These tables contain Monte Carlo analysis where 100 simulations have been done for each target distribution.
Tables report the average ratio of the sup norm (and AMSE) of the IFSs over the corresponding sup norm (respectively AMSE) of the empirical distribution function.

\begin{figure}
\includegraphics{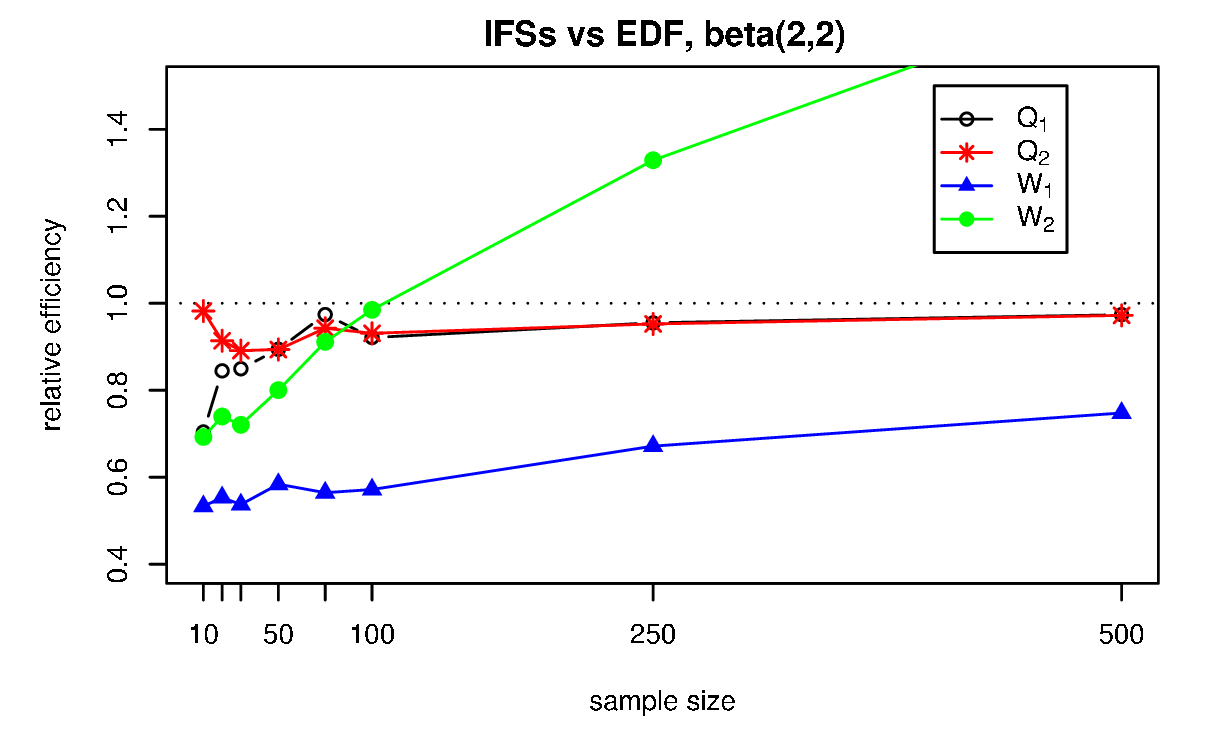}
\includegraphics{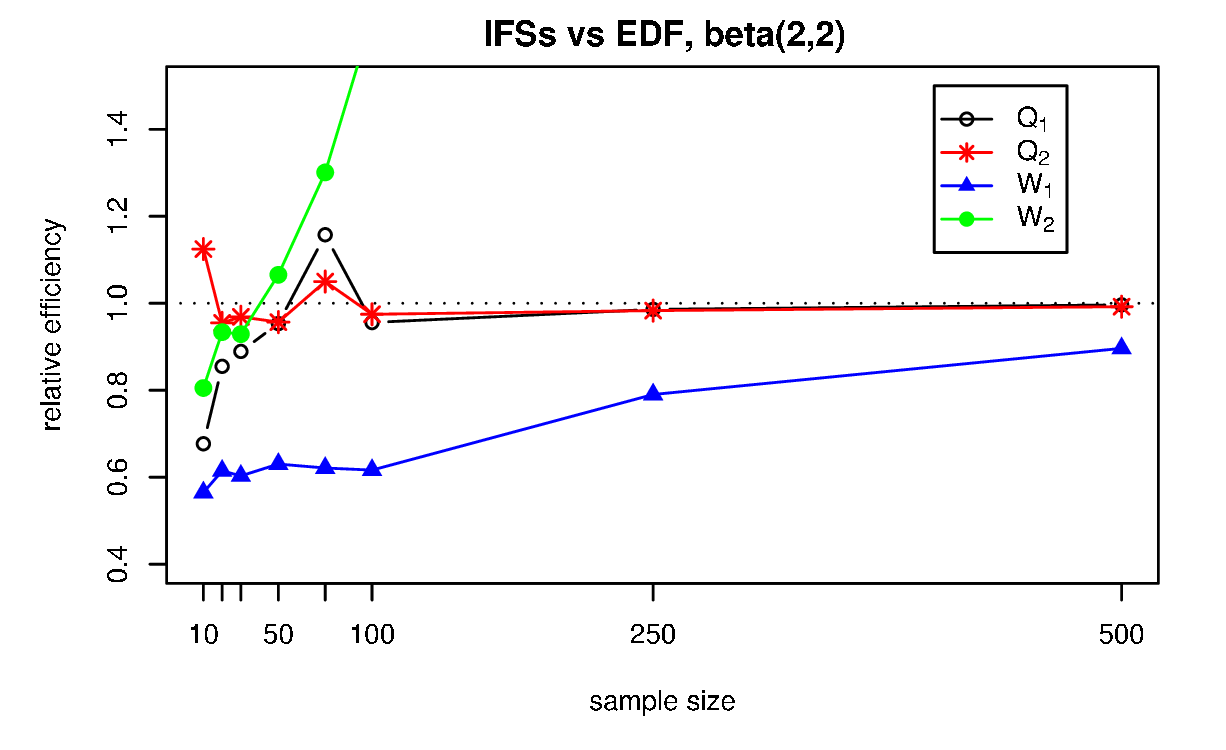}
\caption{Relative efficiency of IFS estimator for different set of maps ${\mathcal W_1}$, ${\mathcal W_2}$,  ${\mathcal Q_1}$ and ${\mathcal Q_2}$ with respect to the empirical distribution function. Based on 100 Monte Carlo simulations. SUP-norm up, AMSE bottom.}
\label{fig:beta22}
\end{figure}

It is possible to notice that the IFS estimator based on maps $\mathcal W_1$ has good properties for  symmetric bell-shaped distributions and distributions with not so heavy tails (see also Figure \ref{fig:beta22}).
It is also evident the asymptotic equivalence of the IFSs to the e.d.f. when quantile maps are used.
Remark that, for $\mathcal W_1$ we have decided to use 62 maps, for $\mathcal W_2$ 28 maps and $n/2$ quantiles for the quantiles maps $\mathcal Q_1$ and $\mathcal Q_2$.
So it is evident that for wavelet-type maps an adjustment can be done by choosing a suitable number of maps in terms of the sample size $n$.

\subsection{What if data are missing?}
Suppose now that the for some reason, the $n$ sample observations from $F$ are in fact a subset of a biggest sample, of unknown size. In practice we do not observe the data on the whole support of $F$ $[\alpha,\beta]$ but only on some windows. This sample reduction has happened due to some sort of censoring. So we are in presence of missing data when we do not know how many data are missing and where exactly they were missed, i.e. we are not in a classical censoring setup. A motivation for this scheme of (non)-observation is the following: suppose one wants to estimate the distribution of the angle of the wind registered by some instruments in degrees (0,360). For some reason, data from angles (15,37) and (62,79) are missing for technical failures or physical obstacles. In this case the empirical distribution function will be flat on these windows and a kernel density estimator will probably show a bimodal behaviour.

\begin{figure}
\includegraphics{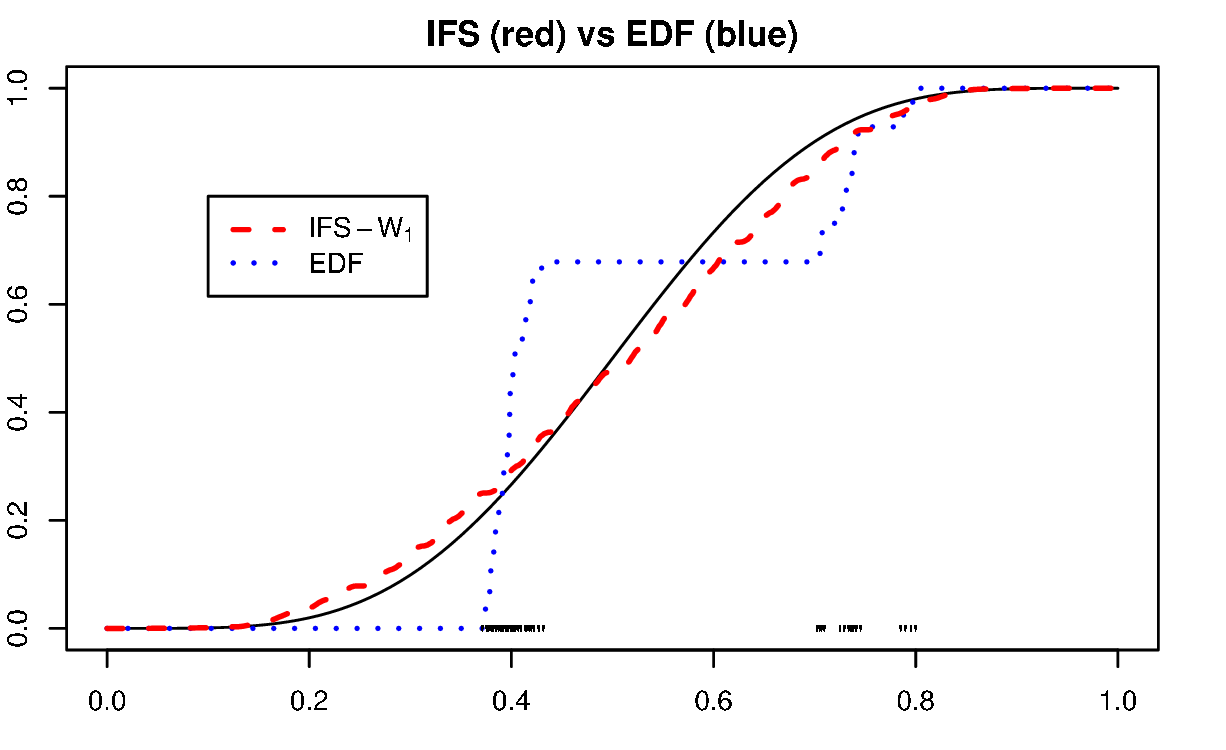}
\includegraphics{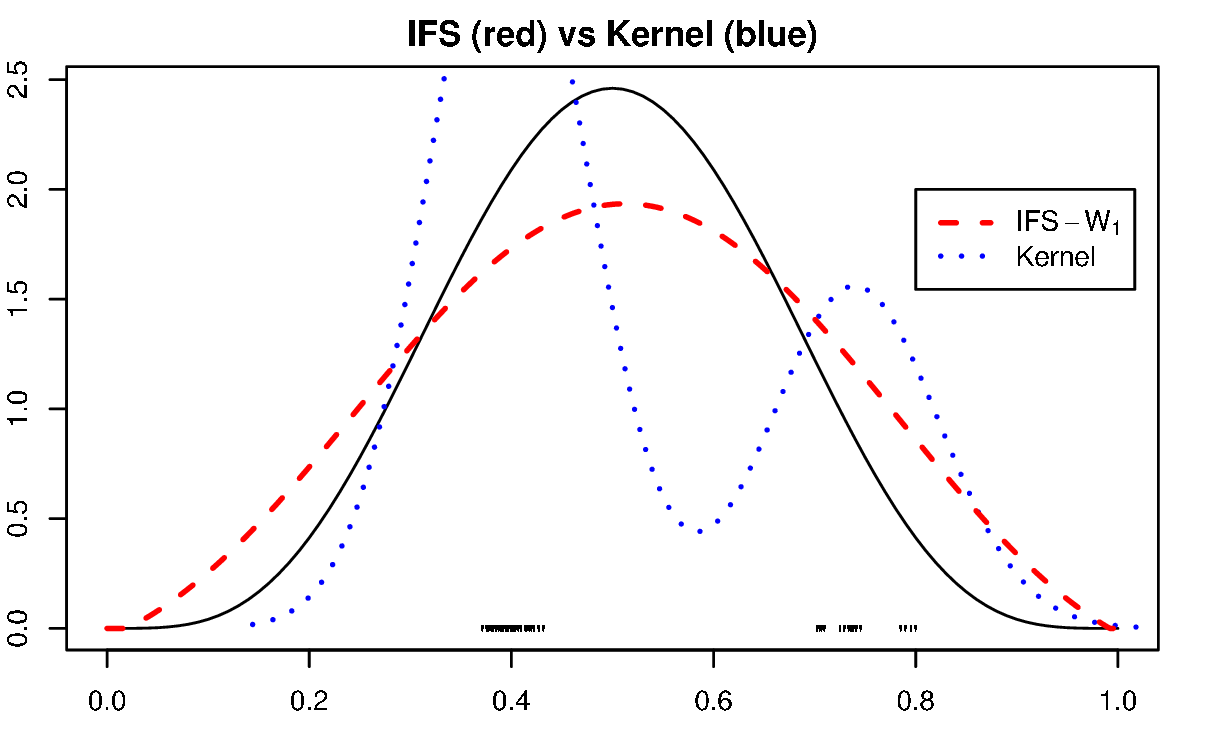}
\caption{Data from a Beta(2,2) distribution when only the observation in $(.1, .15)\cup (.37, .43) \cup (.7, .8)$ are available to the observer all the other being truncated by the instrument.
The observations are marked as vertical ticks. The IFS estimator with $\mathcal W_1$ maps seems to be able to reconstruct the underlying distribution and density function, whistle, for obvious reasons both the edf and the kernel estimators fail. Notice that the arbitrary choice of the window of observation can be changed without substantial loss or gain. In this example the relative efficiency (IFS/EDF) is 7\% for the AMSE and 23\% for the SUP-norm.}
\label{fig:missing}
\end{figure}

Heuristically, this is due to the fact that quantile estimation is inappropriate in this context. At the same time, moments estimation tend to be more robust, in particular if the distribution is symmetric.
We only report a graphical example of what can happen. Figure \ref{fig:missing} is about a sample from a Beta(2,2) distribution when only the observation in $(.1, .15)\cup (.37, .43) \cup (.7, .8)$ are available to the observer all the other being truncated by the instrument (we have choose this interval by hazard). The IFS estimator with $\mathcal W_1$ maps seems to be able to reconstruct the underlying distribution and density function, whistle, for obvious reasons both the e.d.f. and the kernel estimators fail. In this example the relative efficiency (IFS/EDF) is 7\% for the AMSE and 23\% for the SUP-norm which is dramatically better than expected!

\subsection{Algorithm flow for estimation}
\begin{enumerate}
\item calculate sample moments
\item choose the family of maps $\mathcal W$
\item build the quadratic form and solve it for $p$
\item if you want to estimate $F$ at point $x$: take any distribution function, for example the uniform over $[\alpha,\beta]$ and start to iterate $T$
\item stop after few iteration (normally 5 is enough)
\item the ``fixed point'' of $T$ evaluated in $x$ is the estimate of $F(x)$
\end{enumerate}
In case the support of $F$ is not known one case use the range of the sample but the resulting IFS estimator will then try to approximate a distribution function which has exactly that support.
If any hints on the shape of the distribution $F$ is available, use it to choose the maps.

All the examples, tables and graphics have been done by some software developed by the authors. In particular, a package called \texttt{ifs} is freely available for the R environment system \cite{R} in the CRAN (Comprehensive R Archive Network) \texttt{http://cran.R-project.org} under the \textsl{contributed} section.

\section*{Conclusions}
It seems that this kind of approach can be used  to make nonparametric inference when data are missing or sample size are small.
Remark that with this method  it is only possible to work with distributions with compact support. Moreover, a knowledge on the support itself it is needed.
Neverthless, it seams a promising approach and the use of different sets of maps merits further investigation.

% \appendix

% \section{}
% \label{}

% Bibliographic references with the natbib package:
% Parenthetical: \cite{Bai92} produces (Bailyn 1992).
% Textual: \cite{Bai95} produces Bailyn et al. (1995).
% An affix and part of a reference:
%   \cite[e.g.][Ch. 2]{Bar76}
%   produces (e.g. Barnes et al. 1976, Ch. 2).

\eject

\begin{table}
{\scriptsize
\begin{tabular}{c c c}
parameters & AMSE & SUP \\
\begin{tabular}{c|c}
$n$ & law\\
\hline
10 & beta(.9,.1)\\
 10    & beta(.1,.9)\\
  10  & beta(.1,.1)\\
 10 & beta(\,2,\,2)\\
10&beta(\,5,\,5)\\
10&beta(\,3,\,5)\\
10&beta(\,5,\,3)\\
10&beta(\,1,\,1)
\end{tabular}&
\begin{tabular}{c|c|c|c}
${\mathcal W_1} $ & ${\mathcal W_2} $ & ${\mathcal Q_1} $ & ${\mathcal Q_2} $\\
\hline
81.08 & 77.05 & 203.53 & 149.68\\
211.78 & 2024.68 & 203.39 & 258.88\\
118.27 & 416.17 & 182.88 & 104.07\\
56.47 & 80.53 & 67.68 & 112.46\\
52.77 &57.90 &110.35 & 152.29\\
55.95 & 71.07 & 99.92 & 142.52\\
52.50 & 57.34 & 91.75 & 131.37\\
73.35 & 119.04 & 79.01 & 102.04\\
\end{tabular}
&
\begin{tabular}{c|c|c|c}
${\mathcal W_1} $ & ${\mathcal W_2} $ & ${\mathcal Q_1} $ & ${\mathcal Q_2} $\\
\hline
85.76 & 75.44 & 110.11 &110.81\\
175.32 & 441.32 & 114.51 & 161.55\\
114.87 & 192.94 & 119.57 & 106.56\\
53.31 & 69.24 & 70.36 & 98.21\\
53.99 & 54.83 & 81.61 & 125.67\\
51.93 & 60.58 & 81.72 & 116.79\\
51.74 & 52.47 & 77.97 & 109.84\\
65.63 & 90.40 & 70.89 & 90.85\\
\end{tabular}
\end{tabular}
\par
\vspace{12pt}
\par
\begin{tabular}{c c c}
parameters & AMSE & SUP \\
\begin{tabular}{c|c}
$n$ & law\\
\hline
20 & beta(.9,.1)\\
 20    & beta(.1,.9)\\
  20  & beta(.1,.1)\\
 20 & beta(\,2,\,2)\\
20&beta(\,5,\,5)\\
20&beta(\,3,\,5)\\
20&beta(\,5,\,3)\\
20&beta(\,1,\,1)
\end{tabular}&
\begin{tabular}{c|c|c|c}
${\mathcal W_1} $ & ${\mathcal W_2} $ & ${\mathcal Q_1} $ & ${\mathcal Q_2} $\\
\hline
94.69 &85.25 &201.85 &169.78\\
388.83 & 4183.36 & 203.70 & 195.36\\
154.1 & 690.08 & 125.35 & 97.53\\
61.46 & 93.37 &  85.46 &  95.49\\
54.31 & 52.89 & 105.84 & 131.84\\
60.42 & 67.33 & 93.30 &118.51\\
53.82 & 57.72 & 92.26 & 114.84\\
95.93 & 89.79 & 71.66 &  154.54\\ 
\end{tabular}
&
\begin{tabular}{c|c|c|c}
${\mathcal W_1} $ & ${\mathcal W_2} $ & ${\mathcal Q_1} $ & ${\mathcal Q_2} $\\
\hline
90.30 & 79.92 & 105.02 & 123.28\\
257.13 & 612.55 & 109.10 & 122.99\\
139.65 & 255.26 & 103.56 & 99.28\\
55.34 & 73.95 & 84.42 & 91.38\\
53.76 & 48.73 & 85.85 & 106.27\\
55.98 & 60.88 & 85.39 & 101.16\\
53.46 & 52.20 & 85.23 & 102.85\\
63.20 & 106.95 & 81.56 & 82.54\\
\end{tabular}
\end{tabular}
\par
\vspace{12pt}
\par
\begin{tabular}{c c c}
parameters & AMSE & SUP \\
\begin{tabular}{c|c}
$n$ & law\\
\hline
30 & beta(.9,.1)\\
 30    & beta(.1,.9)\\
  30  & beta(.1,.1)\\
 30 & beta(\,2,\,2)\\
30&beta(\,5,\,5)\\
30&beta(\,3,\,5)\\
30&beta(\,5,\,3)\\
30&beta(\,1,\,1)\\
\end{tabular}&
\begin{tabular}{c|c|c|c}
${\mathcal W_1} $ & ${\mathcal W_2} $ & ${\mathcal Q_1} $ & ${\mathcal Q_2} $\\
\hline
107.46 & 90.27 & 195.39 & 143.00\\
540.73 & 6462.03 & 190.82 & 213.45\\
112.66 & 97.04 &  233.50 & 1342.44\\
60.30 & 92.92 & 88.90 & 96.88\\
62.04 & 56.07 & 100.26 & 121.41\\
70.31 & 76.90 & 93.02 & 108.76\\
55.78 & 56.85 & 92.10 & 102.02\\
71.88 & 211.28 & 94.36 &  88.17\\
\end{tabular}
&
\begin{tabular}{c|c|c|c}
${\mathcal W_1} $ & ${\mathcal W_2} $ & ${\mathcal Q_1} $ & ${\mathcal Q_2} $\\
\hline
101.83 & 81.05 & 108.59 & 109.85\\
107.80 & 137.26 & 314.53 & 759.57\\
186.70 & 356.91 & 103.39 & 99.98\\
53.71 & 72.06 &  84.92 & 89.11\\
60.08 & 51.82 &89.26 & 100.16\\
61.68 & 66.29 & 86.36 & 95.24\\
55.56& 51.21 & 88.20 & 94.75\\
63.15 & 121.23 & 83.74 & 83.40\\
\end{tabular}
\end{tabular}
}
\caption{Relative efficiency of IFS estimators with different set of maps  ${\mathcal W_1}$, ${\mathcal W_2}$,  ${\mathcal Q_1}$ and ${\mathcal Q_2}$ with respect to the empirical distribution function (i.e. IFS/EDF).  Based on 100 Monte Carlo simulation for each distribution. Small sample sizes.}
\label{tab:a}
\end{table}

\begin{table}
{\scriptsize
\begin{tabular}{c c c}
parameters & AMSE & SUP \\
\begin{tabular}{c|c}
$n$ & law\\
\hline
50 & beta(.9,.1)\\
 50    & beta(.1,.9)\\
  50  & beta(.1,.1)\\
 50 & beta(\,2,\,2)\\
50&beta(\,5,\,5)\\
50&beta(\,3,\,5)\\
50&beta(\,5,\,3)\\
50&beta(\,1,\,1)
\end{tabular}&
\begin{tabular}{c|c|c|c}
${\mathcal W_1} $ & ${\mathcal W_2} $ & ${\mathcal Q_1} $ & ${\mathcal Q_2} $\\
\hline
132.67 & 115.10 & 163.33 & 129.24\\
1044.12 & 12573.16 & 181.99 & 180.42\\
306.49 & 1917.23 & 105.68 & 97.27\\
63.03 & 106.56 & 95.35 & 95.66\\
68.94 & 60.19 & 102.22 & 114.92\\
79.98 & 93.80 & 96.20 & 102.32\\
63.13 & 62.21 & 93.59 & 98.47\\
73.47 & 304.41 & 97.24 & 92.19\\
\end{tabular}
&
\begin{tabular}{c|c|c|c}
${\mathcal W_1} $ & ${\mathcal W_2} $ & ${\mathcal Q_1} $ & ${\mathcal Q_2} $\\
\hline
109.18& 88.77 & 103.37 & 101.90\\
421.49 & 991.33 & 104.37 & 123.39\\
214.27 & 430.63 & 100.13 &  98.04\\
58.39 & 80.00 & 89.42 & 89.36\\
66.77 & 55.49 & 91.86 & 97.40\\
66.76 & 77.57 & 91.39 & 93.76\\
62.04 & 55.95 & 90.66 & 93.19\\
62.69 & 150.39 & 87.38 & 86.30\\
\end{tabular}
\end{tabular}
\par
\vspace{12pt}
\par
\begin{tabular}{c c c}
parameters & AMSE & SUP \\
\begin{tabular}{c|c}
$n$ & law\\
\hline
100 & beta(.9,.1)\\
 100    & beta(.1,.9)\\
  100  & beta(.1,.1)\\
 100 & beta(\,2,\,2)\\
100&beta(\,5,\,5)\\
100&beta(\,3,\,5)\\
100&beta(\,5,\,3)\\
100&beta(\,1,\,1)
\end{tabular}&
\begin{tabular}{c|c|c|c}
${\mathcal W_1} $ & ${\mathcal W_2} $ & ${\mathcal Q_1} $ & ${\mathcal Q_2} $\\
\hline
195.54 & 158.80 & 140.55 & 108.27\\
1557.30 & 20324.60 & 135.45 & 125.94\\
554.11 & 3918.62 & 102.67 & 98.29\\
61.63 & 165.60 & 95.58 & 97.46\\
87.97 & 67.79 & 99.28 & 108.21\\
111.30 & 134.54&100.68&103.31\\
61.03 & 57.19 & 97.28 & 101.32\\
67.91 & 558.50 & 97.71 & 94.87\\
 \end{tabular}
&
\begin{tabular}{c|c|c|c}
${\mathcal W_1} $ & ${\mathcal W_2} $ & ${\mathcal Q_1} $ & ${\mathcal Q_2} $\\
\hline
138.93 & 105.31 & 102.25 & 99.07\\
536.84 & 1267.81 & 103.87 & 106.05\\
304.59 & 625.75 & 99.10 & 98.04\\
57.18 & 98.50 & 92.11 & 93.09\\
78.94 &60.96 &94.83 &96.52\\
79.59 & 100.20 & 95.35 & 95.72\\
65.97 & 55.08 & 94.14 & 95.42\\
58.71 & 201.10 & 90.83 & 89.97\\
 \end{tabular}
\end{tabular}
\par
\vspace{12pt}
\par
\begin{tabular}{c c c}
parameters & AMSE & SUP \\
\begin{tabular}{c|c}
$n$ & law\\
\hline
250 & beta(.9,.1)\\
 250    & beta(.1,.9)\\
  250  & beta(.1,.1)\\
 250 & beta(\,2,\,2)\\
250&beta(\,5,\,5)\\
250&beta(\,3,\,5)\\
250&beta(\,5,\,3)\\
250&beta(\,1,\,1)\\
\end{tabular}&
\begin{tabular}{c|c|c|c}
${\mathcal W_1} $ & ${\mathcal W_2} $ & ${\mathcal Q_1} $ & ${\mathcal Q_2} $\\
\hline
338.72 & 255.23 & 115.25 & 101.55\\
3979.61 & 50448.13 &117.81 & 105.37\\
1345.72 & 10051.20 &  100.60 & 98.97\\
79.01 & 275.93 & 98.59 & 98.30\\
163.68 & 99.35 &99.07 & 100.54\\
212.17 & 228.58 & 99.45 & 99.69\\
91.32 & 73.31 & 99.05 &99.20\\
69.03 & 1165.61 & 99.47 & 98.46\\ 
 \end{tabular}
&
\begin{tabular}{c|c|c|c}
${\mathcal W_1} $ & ${\mathcal W_2} $ & ${\mathcal Q_1} $ & ${\mathcal Q_2} $\\
\hline
180.29 & 131.97 &  100.68 &  99.43\\
874.65 & 2045.15 & 100.82 & 99.73\\
480.12 & 977.30 & 99.16 & 98.73\\
67.14 & 132.87 & 95.50 & 95.24\\
111.38 & 78.48 & 96.40 & 96.83\\
113.70 & 142.21 & 96.57 & 96.32\\
88.87 & 67.13 & 96.84 & 97.24\\
61.07 & 293.58 & 94.88 & 94.55\\
 \end{tabular}
\end{tabular}
}
\caption{Relative efficiency of IFS estimators with different set of maps  ${\mathcal W_1}$, ${\mathcal W_2}$,  ${\mathcal Q_1}$ and ${\mathcal Q_2}$ with respect to the empirical distribution function  (i.e. IFS/EDF). Based on 100 Monte Carlo simulation for each distribution. Moderate to big sample sizes.}
\label{tab:b}
\end{table}

\begin{figure}
\includegraphics{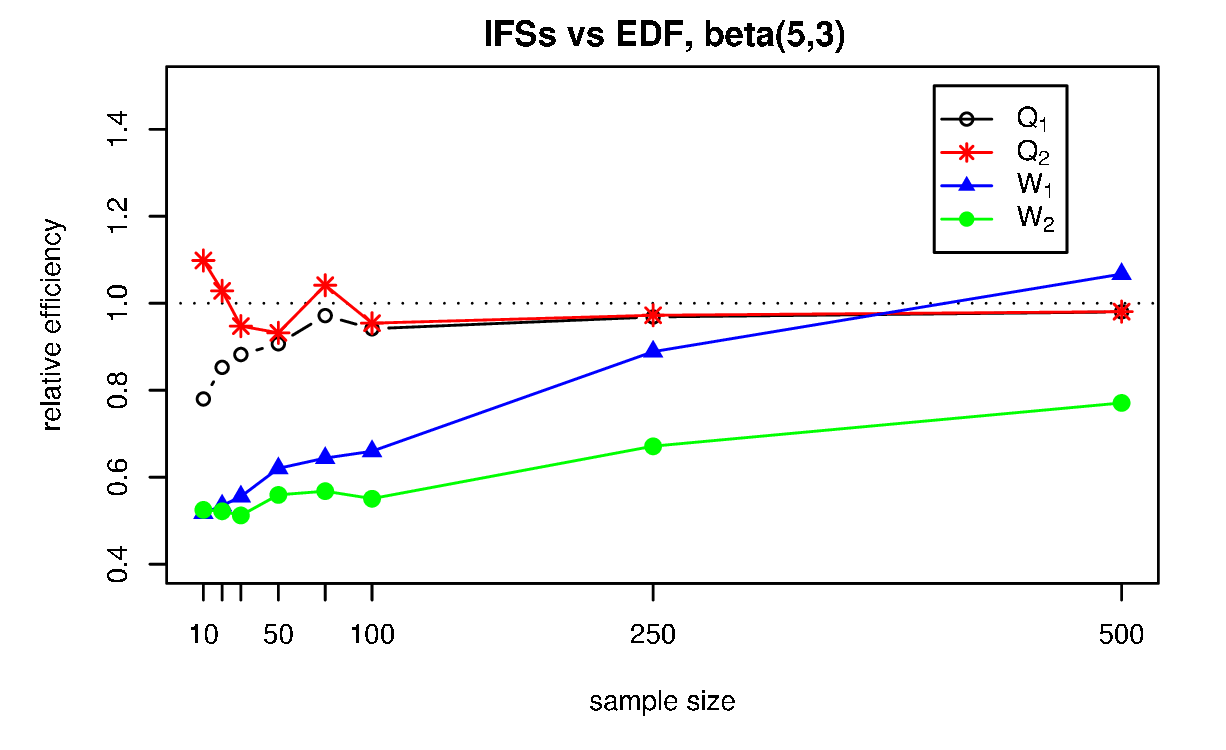}
\includegraphics{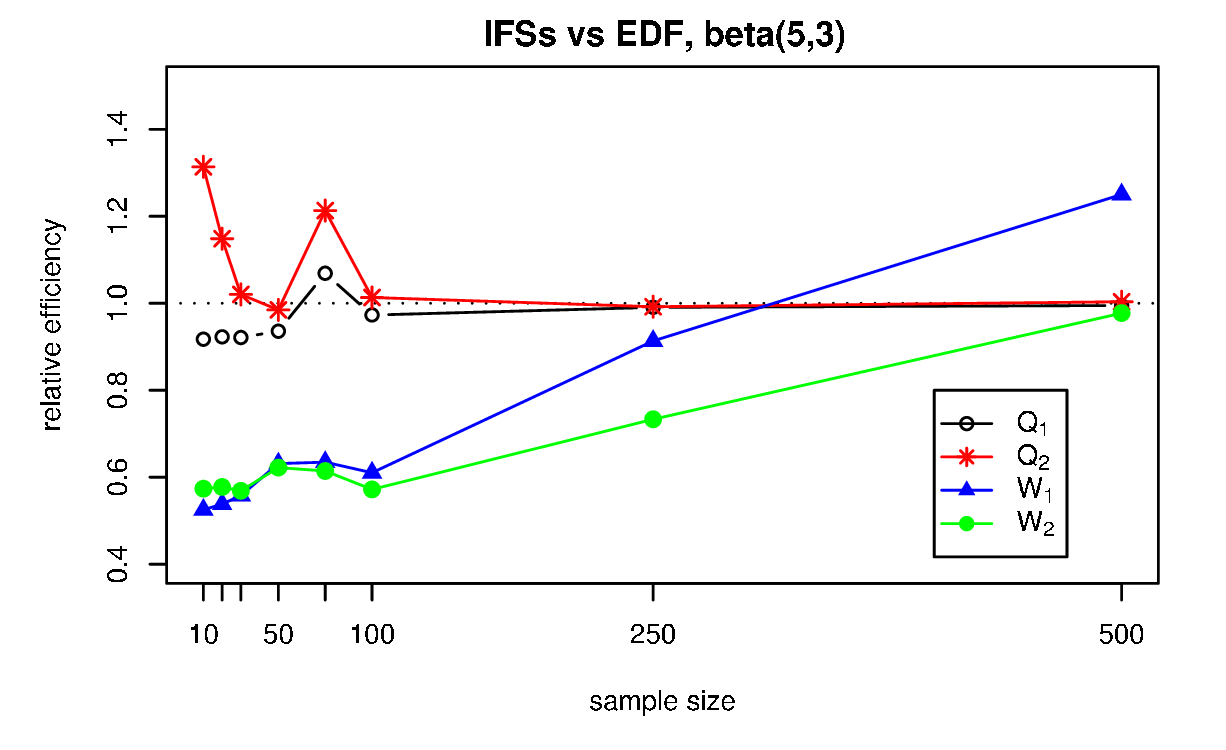}
\caption{Relative efficiency of IFS estimator for different set of maps ${\mathcal W_1}$, ${\mathcal W_2}$,  ${\mathcal Q_1}$ and ${\mathcal Q_2}$ with respect to the empirical distribution function. Based on 100 Monte Carlo simulations. SUP-norm up, AMSE bottom.}
\label{fig:beta53}
\end{figure}

\begin{figure}
\includegraphics{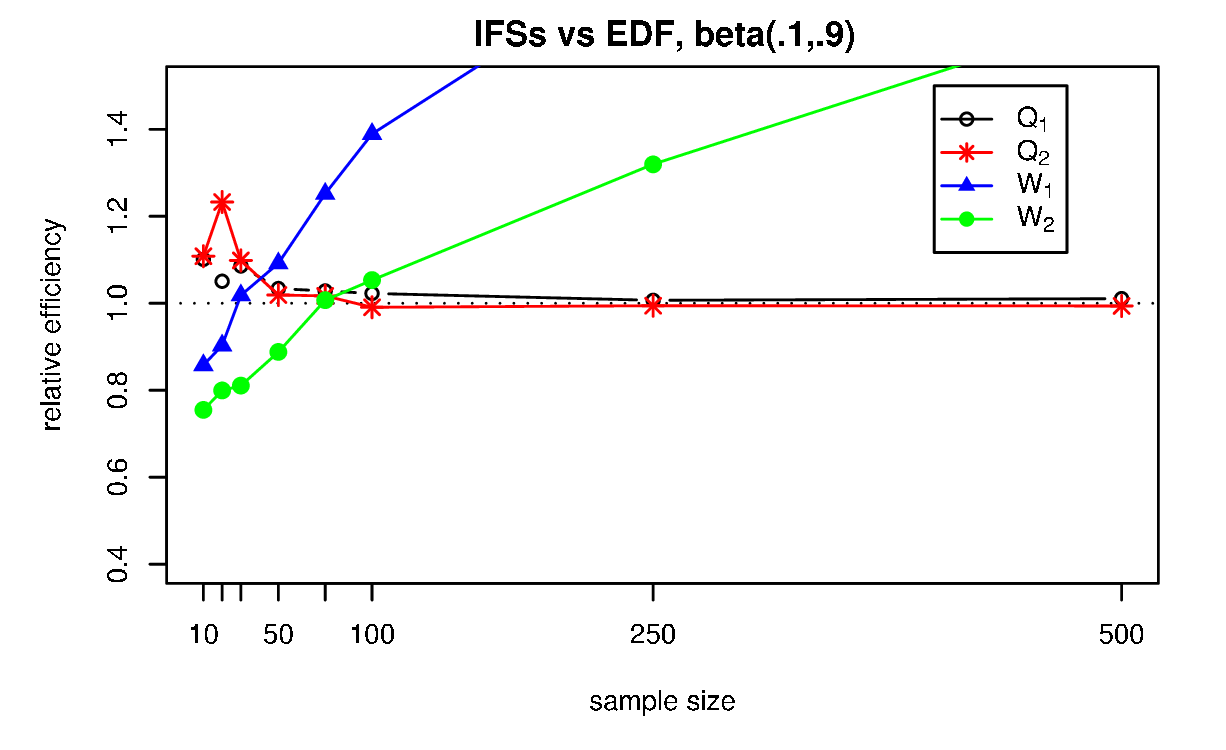}
\includegraphics{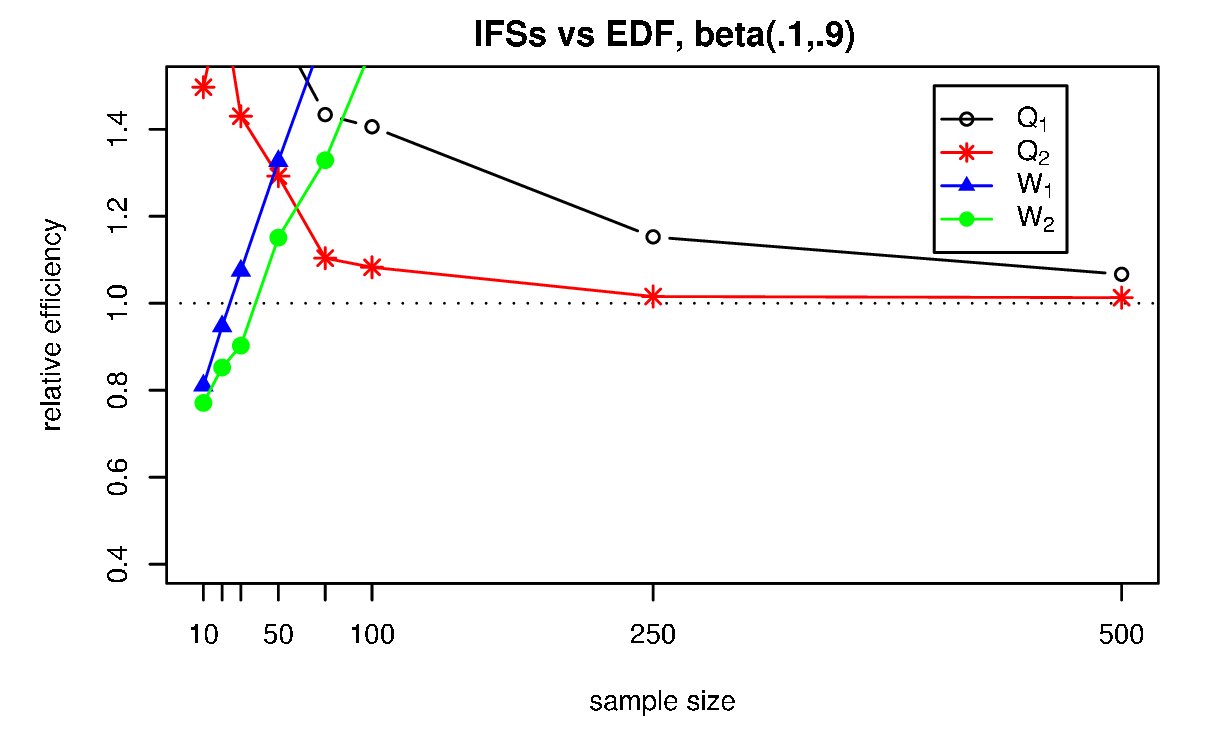}
\caption{Relative efficiency of IFS estimator for different set of maps ${\mathcal W_1}$, ${\mathcal W_2}$,  ${\mathcal Q_1}$ and ${\mathcal Q_2}$ with respect to the empirical distribution function. Based on 100 Monte Carlo simulations. SUP-norm up, AMSE bottom.}
\label{fig:beta19}
\end{figure}

\end{document}